\documentclass{amsart}
\usepackage{amsmath,amsfonts,latexsym,amssymb}

\newtheorem{theorem}{Theorem}
\newtheorem{lemma}{Lemma}

\newtheorem{proposition}{Proposition}

\theoremstyle{remark}

\newtheorem*{ack}{Acknowledgements}

\begin{document}

\markboth{Ritabrata Munshi}{$t$-aspect subconvexity for $GL(3)$ $L$-functions}
\title[$t$-aspect subconvexity for $GL(3)$ $L$-functions]{The circle method and bounds for $L$-functions - III:\\ $t$-aspect subconvexity for $GL(3)$ $L$-functions}

\author{Ritabrata Munshi}   
\address{School of Mathematics, Tata Institute of Fundamental Research, 1 Dr. Homi Bhabha Road, Colaba, Mumbai 400005, India.}     
\email{rmunshi@math.tifr.res.in}


\begin{abstract}
Let $\pi$ be a Hecke-Maass cusp form for $SL(3,\mathbb Z)$. In this paper we will prove the following subconvex bound
$$
L\left(\tfrac{1}{2}+it,\pi\right)\ll_{\pi,\varepsilon} (1+|t|)^{\frac{3}{4}-\frac{1}{16}+\varepsilon}.
$$
\end{abstract}


\maketitle


\section{Introduction}
\label{intro}

Let $\pi$ be a Hecke-Maass cusp form of type $(\nu_1,\nu_2)$ for $SL(3,\mathbb Z)$. Let the normalized Fourier coefficients of $\pi$ be given by $\lambda(m_1, m_2)$ (so that $\lambda(1, 1)=1$). The Langlands parameters $(\alpha_1, \alpha_2, \alpha_3)$ associated with $\pi$ are defined as $\alpha_1=-\nu_1-2\nu_2+1$, $\alpha_2=-\nu_1+\nu_2$ and $\alpha_3=2\nu_1+\nu_2-1$. The Ramanujan-Selberg conjecture predicts that $\text{Re}(\alpha_i)=0$. From the work of Jacquet and Shalika \cite{JS}, we (at least) know that $|\text{Re}(\alpha_i)|<\frac{1}{2}$. The $L$-series associated with $\pi$ is given by 
$$
L(s,\pi)=\sum_{n=1}^\infty \lambda(1,n)n^{-s}
$$
in the domain $\sigma=\text{Re}(s)>1$. This extends to an entire function and satisfies a functional equation. More precisely there is an associated gamma factor given by
$$
\gamma(s,\pi)=\prod_{i=1}^3\pi^{-\frac{s}{2}}\Gamma\left(\frac{s-\alpha_i}{2}\right),
$$ 
so that
$$
\gamma(s,\pi)L(s,\pi)=\gamma(s,\tilde\pi)L(1-s,\tilde\pi).
$$ 
Here $\tilde\pi$ is the dual form having Langlands parameters $(-\alpha_3, -\alpha_2, -\alpha_1)$. The convexity principle implies that $L(1/2+it,\pi)\ll_{\pi} (1+|t|)^{3/4}$ - the convexity bound. The purpose of this paper is to prove the following.
\begin{theorem}
\label{mthm}
Let $\pi$ be a Hecke-Maass cusp form for $SL(3,\mathbb Z)$. Then we have
$$
L\left(\tfrac{1}{2}+it,\pi\right)\ll_{\pi,\varepsilon} (1+|t|)^{\frac{3}{4}-\frac{1}{16}+\varepsilon}.
$$
\end{theorem}

A similar subconvex bound, with same exponent, is known for the symmetric square lifts of $SL(2,\mathbb Z)$ forms (or self dual forms for $SL(3,\mathbb Z)$) due to the work of Li \cite{L}. Other subconvexity results in the case of degree three $L$-functions in different aspects can be found in \cite{B}, \cite{Mu1}, \cite{Mu2}, \cite{Mu3} and \cite{Mu0}.  Subconvex bound in the $t$-aspect was first established by Weyl \cite{W} for degree one $L$-functions, and by Good \cite{Go} for degree two $L$-functions. This paper settles the problem for degree three $L$-functions. \\

Like the two previous papers \cite{Mu} and \cite{Mu0}, with the same title, we will yet again demonstrate the power of the circle method in the context of subconvexity. In the present situation Kloosterman's version of the circle method works best. Let
$$
\delta(n)=\begin{cases} 1&\text{if}\;\;n=0;\\
0&\text{otherwise}.\end{cases}
$$
Then for any real number $Q$, we have
\begin{align}
\label{cm}
\delta(n)=2\:\text{Re}\int_0^1\mathop{\sum\sideset{}{^\star}\sum}_{1\leq q\leq Q <a\leq q+Q}\frac{1}{aq}e\left(\frac{n\bar a}{q}-\frac{nx}{aq}\right)\mathrm{d}x
\end{align}
for $n\in\mathbb Z$ (and $e(z)=e^{2\pi iz}$). The $\star$ on the sum indicates that the sum over $a$ is restricted by the condition $(a,q)=1$, also $\bar a$ stands for the multiplicative inverse of $a$ modulo $q$. (For a proof of this formula see \cite{IK}.) There are well understood drawbacks in this form of circle method. However in our treatment these do not create any problem. After an application of the Poisson summation formula, we will be able to write $a$ in terms of the dual frequency (see \eqref{a-det} in Subsection~\ref{poi-1}), and hence we do not need to execute the complete character sum over $a$. After that we will only need the fact that $a\asymp Q$. (The notation $\alpha\asymp A$ means that there exists absolute constants $0<c_1<c_2$ such that $c_1A< |\alpha| <c_2A$.) The main advantage of the above version of the circle method is the explicit form of the weight function $e(-nx/aq)$, which will be helpful for estimating the exponential integrals in Section~\ref{aoi}. \\

Suppose $t>2$, then by approximate functional equation  we have
\begin{align}
\label{afe}
L\left(\tfrac{1}{2}+it,\pi\right)\ll t^{\varepsilon}\; \mathop{\sup}_{N\leq t^{3/2+\varepsilon}} \frac{|S(N)|}{N^{1/2}}+t^{-2012}
\end{align}
where $S(N)$ is a sum of type 
$$
S(N):=\sum_{n=1}^\infty \lambda(1,n)n^{-it}V\left(\frac{n}{N}\right)
$$
for some smooth function $V$ (allowed to depend on $t$) supported in $[1,2]$ and satisfying $V^{(j)}(x)\ll_j 1$. (In this paper the notation $\alpha\ll A$ will mean that there is a constant $c$ such that $|\alpha|\leq cA$. The dependence of the constant on the automorphic form $\pi$, and $\varepsilon$, when occurring, will be ignored.)  Hence to establish subconvexity we need to show cancellation in the sum $S(N)$ for $N$ roughly of size $t^{3/2}$. We can and shall further normalize $V$, for convenience, so that $\int V(y)\mathrm{d}y=1$.\\

We apply \eqref{cm} directly to $S(N)$ as a device for separation of the oscillation of the Fourier coefficients $\lambda(1,n)$ and $n^{-it}$. This by itself does not seem very effective, and as in \cite{Mu0} we need a conductor lowering mechanism. For this purpose we introduce an extra integral namely
$$
S(N)=\frac{1}{K}\int_{\mathbb R}V\left(\frac{v}{K}\right)\mathop{\sum\sum}_{\substack{n,m=1\\n=m}}^\infty \lambda(1,n)m^{-it}\:\left(\frac{n}{m}\right)^{iv}V\left(\frac{n}{N}\right)U\left(\frac{m}{N}\right)\mathrm{d}v,
$$
where $t^\varepsilon<K<t$ is a parameter which will be chosen optimally later, and  $U$ is a smooth function supported in $[1/2,5/2]$, with $U(x)=1$ for $x\in[1,2]$ and $U^{ (j)}\ll_j 1$. For $n,m\in [N,2N]$, the integral
$$
\frac{1}{K}\int_{\mathbb R}V\left(\frac{v}{K}\right)\left(\frac{n}{m}\right)^{iv}\mathrm{d}v
$$
is small if $|n-m|\gg Nt^{\varepsilon}/K$. This step is the most crucial `trick' in this paper. The reader should realize that this is the analytic analogue of the arithmetic condition $M_1|(n-m)$, that we had in \cite{Mu0}. There we used this to replace $\delta(n-m)$ by $\delta((n-m)/M_1)$, and hence the optimal choice of the size of the modulus in the circle method reduced from $N^{1/2}$ (where $n,m\asymp N$) to $(N/M_1)^{1/2}$. Here also it reduces the optimal size of the modulus in the circle method from $N^{1/2}$ to $(N/K)^{1/2}$. This is probably not completely obvious at this stage, but it should become clear later.\\

So the `natural choice' for $Q$ in \eqref{cm} to detect the event $n-m=0$ is 
\begin{align}
\label{q-choice}
Q=\left(\frac{N}{K}\right)^{1/2}
\end{align}
and we get 
$$
S(N)=S^+(N)+S^-(N)
$$
where
\begin{align}
\label{S(N)-cm}
S^\pm(N)=&\frac{1}{K}\int_0^1\int_{\mathbb R}V\left(\frac{v}{K}\right)\mathop{\sum\sideset{}{^\star}\sum}_{1\leq q\leq Q <a\leq q+Q}\frac{1}{aq}\\
\nonumber &\times \mathop{\sum\sum}_{n,m=1}^\infty \lambda(1,n)n^{iv}m^{-i(t+v)}e\left(\pm\frac{(n-m)\bar a}{q}\mp\frac{(n-m)x}{aq}\right)V\left(\frac{n}{N}\right)U\left(\frac{m}{N}\right)\mathrm{d}v\mathrm{d}x.
\end{align} \\

In the rest of the paper we will analyse $S^+(N)$ (the same analysis holds for $S^-(N)$), with $Q$ as in \eqref{q-choice}, using summation formulae and stationary phase method. We will take 
\begin{align}
\label{cond-K}
t^{11/8}<N<t^{3/2+\varepsilon},\;\;\;\;\text{and}\;\;\;\; \frac{t^{6/5}}{N^{3/5}}\leq K<\min\left\{\frac{t^{2-\varepsilon}}{N}, N^{1/3}\right\}.
\end{align}
The optimal choice of $K$, as we will see at the end, is given by $K=\max\{N^{1/4},t^{6/5}/N^{3/5}\}$. With this choice of $K$ we will establish the following bound.\\

\begin{proposition}
\label{prop1}
We have
\begin{align}
\label{prop-eq}
S^+(N)\ll \begin{cases} t^{1/2+\varepsilon} N^{5/8} &\text{if}\;\;\;t^{24/17}<N\ll t^{3/2+\varepsilon};\\
t^{11/10+\varepsilon}N^{1/5} & \text{if}\;\;\;t^{11/8}<N\leq t^{24/17}.\end{cases}
\end{align}\\
\end{proposition}

Same bound holds for $S^-(N)$, and consequently, for $S(N)$. For $N\leq t^{11/8}$ the trivial bound $S(N)\ll Nt^\varepsilon$, which follows from Lemma \ref{ram-on-av} (i.e. Ramanujan bound on average) of Section~\ref{prelim}, is sufficient for our purpose. Clearly Theorem \ref{mthm} follows from \eqref{afe} and \eqref{prop-eq} (after a short computation). In the rest of the paper we will prove the proposition. 
\\

Let us now briefly explain the steps in the proof. Temporarily assume the Ramanujan conjecture $\lambda(1,n)\ll n^\varepsilon$. This is not very serious, as at any step where it is required one can use Cauchy inequality and use Ramanujan bound on average, i.e. Lemma~\ref{ram-on-av}. The circle method has been used to separate the sums on $n$ and $m$ and we have arrived at \eqref{S(N)-cm}. Trivially estimating the sum we get $S(N)\ll N^{2+\varepsilon}$. (For simplicity assume that $N=t^{3/2}$ and $q\asymp Q$.) So we are required to save $N$ (and a little more) in a sum of the form 
\begin{align*}
\int_{K}^{2K} \sum_{q\asymp Q}\mathop{\sideset{}{^\star}\sum}_{Q <a\leq q+Q} \mathop{\sum}_{n\asymp N} \lambda(1,n)n^{iv}e\left(\frac{n\bar a}{q}-\frac{nx}{aq}\right)\sum_{m\asymp N}m^{-i(t+v)}e\left(-\frac{m\bar a}{q}+\frac{mx}{aq}\right)\mathrm{d}v.
\end{align*}
The sum over $m$ has `conductor' $Qt\asymp N^{1/2}t/K^{1/2}$. Roughly speaking the conductor takes into account both the arithmetic modulus, which is $q$, and the amplitude of oscillation in the analytic weight function, which is of size $t$. Note that both $m^{-i(t+v)}$ and $m^{-it}$ have same amplitude if $|v|\ll t^{1-\varepsilon}$. So the extra oscillating term, namely $m^{-iv}$ which we are inserting is not hurting us here. On the other hand larger is the $K$ smaller is the arithmetic modulus. So the overall conductor in the sum over $m$ is reduced. Applying Poisson summation (and `executing' the sum over $a$) we are able to save $N/(Qt)^{1/2}\times Q^{1/2}=N/t^{1/2}$. Of course to this end we need the second derivative bound for the resulting exponential integral. In fact we need to use the stationary phase method. Observe that the saving so far is independent of $K$. Now we need to save $t^{1/2}$ in a sum of the form 
\begin{align*}
\int_{K}^{2K}\mathop{\sum}_{q\asymp Q}\mathop{\sum}_{\substack{(m,q)=1\\|m|\asymp Qt/N}}\:\left(\frac{(t+v)aq}{(x-ma)}\right)^{-i(t+v)}\:\mathop{\sum}_{n\asymp N} \lambda(1,n)e\left(\frac{nm}{q}\right)n^{iv}e\left(-\frac{nx}{aq}\right)\mathrm{d}v,
\end{align*}
where $a$ is the unique multiplicative inverse of $m$ modulo $q$ in the range $(Q,q+Q]$.\\

Consider the sum over $n$, which involves the Fourier coefficients, and has `conductor' $(QK)^3$. Observe that larger values of $K$ is taking us to a worse situation. But applying Voronoi summation formula we are able to save $N/(QK)^{3/2}=N^{1/4}/K^{3/4}$. To this end we need Weil bound for Kloosterman sums and second derivative bound for certain exponential integrals that arise in the integral transform resulting from Voronoi. Moreover we are able to save $K^{1/2}$ in the integral over $v$ (see Section~\ref{aoi}). We now need to save $K^{1/4}t^{1/8}$ in a sum of the form 
\begin{align*}
\mathop{\sum}_{n\asymp N^{1/2}K^{3/2}} \lambda(n,1)\mathop{\sum\sum}_{\substack{q\asymp Q,\; (m,q)=1\\|m|\asymp qt/N}}S(\bar m, n; q)\:\int_{-K}^K n^{-i\tau}g(q,m,\tau)\mathrm{d}\tau
\end{align*}
where the function $g$ is of size $O(1)$ but highly oscillatory. It is basically a product of a smooth bump function and 
\begin{align*}
\left(\frac{N}{q^3}\right)^{-i\tau}\left(-\frac{(t+\tau)q}{Nm}\right)^{-i(t+\tau)}.
\end{align*}
The role of the $v$ integral and the parameter $K$ is not yet clear. At this moment it seems to be hurting us more rather than helping. \\

The next step involves taking Cauchy to get rid of the Fourier coefficients, but this process also squares the amount we need to save. So now we face with the task of saving $K^{1/2}t^{1/4}$ in a sum which roughly looks like
$$
\sum_{n\asymp N^{1/2}K^{3/2}}\Bigl|\mathop{\sum\sum}_{\substack{q\asymp Q\\|m|\asymp Qt/N\\ (m,q)=1}}S(\bar m, n;q)\int_{-K}^K n^{-i\tau}g(q,m,\tau)\mathrm{d}\tau\Bigr|^2.
$$
One should note that we need to save $K^{1/2}t^{1/4}$ together with square root saving in the Kloosterman sum and $K^{1/2}$ saving in the integral (which is the second derivative bound). The idea is to open the absolute square and execute the sum over $n$ using Poisson summation. The resulting diagonal contribution or the zero frequency contribution is satisfactory for our purpose if 
$$
\text{the number of terms inside the absolute value}\;= \frac{Q^2t}{N} > K^{1/2}t^{1/4}
$$
or equivalently $t^{1/2}>K$. On the other hand by Poisson we make a saving (ignoring the zero frequency) of size $N^{1/2}K^{3/2}/QK^{1/2}=K^{3/2}$, as the length of the sum is $N^{1/2}K^{3/2}$ and the conductor is of size $Q^2K=N$. This is where we are getting help from the parameter $K$. The conductor is independent of $K$, but the length of the sum increases with $K$. So effectively we have a drop in the conductor of the sum. So the contribution of the non-zero frequencies is satisfactory if
$$
K^{3/2}>K^{1/2}t^{1/4}
$$
or $K>t^{1/4}$. In particular by choosing $K$ in the range $t^{1/4}<K<t^{1/2}$ we can get a bound which breaks the convexity barrier.\\

We conclude this section with a heuristic argument.
We want to bound 
$$
S=\sum_{n\asymp t^{3/2}}\lambda(1,n)n^{it}.
$$
Consider the function $n\mapsto \chi(n)n^{iv}$, where $\chi$ varies over all characters of conductor $\leq Q$ and $v\asymp K$. These functions (if one discretizes $v$ naturally) span a vector space of dimension $\approx Q^2K=t^{3/2}$. Hence one expects to write the function $n\mapsto n^{it}$ as a linear combination of the above set of functions. A priori it is not clear how to actually do so, but the formula \eqref{SNC} of Lemma~\ref{pre-lemma} amounts to a formula of this general nature. We have thus written $S$ in terms of 
$$
\sum_{n\asymp t^{3/2}}\lambda(1,n)\chi(n)n^{iv}.
$$ 
Now we apply the functional equation for the $L$-function $L(1/2-iv,\pi\otimes \chi)$. The sum gets dualized to a sum of length $(QK)^3/t^{3/2}$ (which is the range of summation in Lemma~\ref{lemma-snc}). It is of course essential, for this strategy to work, that this dual length is shorter than the original length, i.e.
$$
(QK)^3/t^{3/2}<t^{3/2}.
$$
Since $Q=t^{3/4}/K^{1/2}$ the inequality boils down to the requirement $K<t^{1/2}$. One should realize that we are trading an $L$-function of conductor $t^3$ for a family of $L$-functions of conductor $(QK)^3$, and this is beneficial only when the above inequality holds. It is also worth observing that - for this exact reason - degree $3$ $L$-function is the limit of the method (assuming we restrict to twisting only by $GL(1)$). Indeed for $L$-function of degree $r$ one will need $K<t^{2-r/2}$, which is not possible for $r\geq 4$.\\

\ack
The author is partly supported by SwarnaJayanti Fellowship, Department of Science and Technology, Government of India. 
He also gratefully acknowledges the hospitality of TIFR CAM Bangalore, where a part of the work was done.\\


\section{$GL(3)$ Voronoi summation formula and stationary phase method}
\label{prelim}

\subsection{Voronoi type summation formula for $SL(3,\mathbb Z)$}
Suppose $\pi$ is a Maass form of type $(\nu_1,\nu_2)$ for $SL_3(\mathbb Z)$, which is an eigenfunction of all the Hecke operators with Fourier coefficients $\lambda(n_1,n_2)$, normalized so that $\lambda(1,1)=1$. Since we shall work entirely at the level of $L$-functions, we simply refer to Goldfeld's book \cite{G} for details regarding automorphic forms on higher rank groups. In this subsection we recall two important results - a summation formula for the Fourier coefficients twisted by additive characters and a bound on the average size of the Fourier coefficients - which will play vital role in our analysis. \\

Let $g$ be a compactly supported smooth function on $(0,\infty)$, and let 
$
\tilde g(s)=\int_0^\infty g(x)x^{s-1}dx
$ 
be its Mellin transform. For $\sigma>-1+\max\{-\text{Re}(\alpha_1),-\text{Re}(\alpha_2),-\text{Re}(\alpha_3)\}$ and $\ell=0,1$ define
$$
\gamma_\ell(s):=\frac{\pi^{-3s-\frac{3}{2}}}{2}\:\frac{\Gamma\left(\frac{1+s+\alpha_1+\ell}{2}\right)\Gamma\left(\frac{1+s+\alpha_2+\ell}{2}\right)\Gamma\left(\frac{1+s+\alpha_3+\ell}{2}\right)}{\Gamma\left(\frac{-s-\alpha_1+\ell}{2}\right)\Gamma\left(\frac{-s-\alpha_2+\ell}{2}\right)\Gamma\left(\frac{-s-\alpha_3+\ell}{2}\right)},
$$
set $\gamma_\pm(s)=\gamma_0(s)\mp i\gamma_1(s)$ and let
\begin{align}
\label{gl}
G_{\pm}(y)=\frac{1}{2\pi i}\int_{(\sigma)}y^{-s}\gamma_\pm(s)\tilde g(-s)ds.
\end{align}
The following Voronoi type summation formula (see \cite{L}, \cite{MS}) will play a crucial role in our analysis. Recall the definition of the Kloosterman sum -
$$
S(a,b; c)=\sideset{}{^\star}\sum_{\alpha\bmod{c}}e\left(\frac{a\alpha+b\overline{\alpha}}{c}\right),
$$ 
where $\bar{\alpha}$ denotes the multiplicative inverse of $\alpha\bmod{c}$.

\begin{lemma}
\label{lem-voronoi}
Let $g$ be a compactly supported smooth function on $(0,\infty)$, we have
\begin{align}
\label{voronoi3}
\sum_{n=1}^\infty \lambda(1,n)e\left(\frac{an}{q}\right)g(n)=&q\sum_\pm\sum_{n_1|q}\sum_{n_2=1}^\infty \frac{\lambda(n_2,n_1)}{n_1n_2}S(\bar a, \pm n_2; q/n_1)G_\pm\left(\frac{n_1^2n_2}{q^3}\right),
\end{align}
where $(a,q)=1$ and $\bar{a}$ denotes the multiplicative inverse of $a\bmod{q}$.\\
\end{lemma}

We need to study the behaviour of the gamma factor $\gamma_\pm(s)$ more closely, especially for $s$ restricted in vertical strips. Using Stirling formula we can pull out the oscillatory part, and the remaining part satisfies a `scaling property'. Indeed for $s=-\frac{1}{2}+i\tau$ with $|\tau|\gg t^\varepsilon$, we apply Stirling's formula to get
\begin{align}
\label{scaling}
\gamma_\pm\left(-\frac{1}{2}+i\tau\right)=\left(\frac{|\tau|}{e\pi}\right)^{3i\tau}\Phi_\pm(\tau),\;\;\;\;\text{where}\;\;\;\;\Phi_\pm'\left(\tau\right)\ll \frac{1}{|\tau|}.
\end{align}\\

The following lemma, which gives Ramanujan conjecture on average, is also well-known. It follows from standard properties of the Rankin-Selberg $L$-function.
\begin{lemma}
\label{ram-on-av}
We have
$$
\mathop{\sum\sum}_{n_1^2n_2\leq x}|\lambda(n_1,n_2)|^2\ll x^{1+\varepsilon},
$$
where the implied constant depends on the form $\pi$ and $\varepsilon$.\\
\end{lemma}

\subsection{Stationary phase method}
We will need estimates for exponential integrals of the form 
\begin{align}
\label{exp-integral}
\mathfrak{I}=\int_{a}^b g(x)e(f(x))\mathrm{d}x,
\end{align}
where $f$ and $g$ are smooth real valued functions. First we recall an easy estimate. Suppose in the range of the integral we have $|f'(x)| \geq B$ and $f^{(j)}(x)\ll B^{1+\varepsilon}$ for $j\geq 2$. Suppose $\text{Supp}(g)\subset (a,b)$ and $g^{(j)}(x)\ll_{a,b,j} 1$. Then by making the change of variable $u=f(x)$ we get
$$
\mathfrak{I}=\int_{f(a)}^{f(b)} \frac{g(f^{-1}(u))}{f'(f^{-1}(u))}e(u)\mathrm{d}u.
$$
By applying integration by parts $j$ times we get
$$
\mathfrak{I}\ll_{a,b,j,\varepsilon} B^{-j+\varepsilon}.
$$
This will be used at several places to show that certain exponential integrals are negligibly small in the absence of the stationary phase.\\

In case there is a single stationary phase point then the integral has an asymptotic expansion. A sharp version of this stationary phase method, which can be found in \cite{Hu}, will be useful for our purpose. The estimates are in terms of parameters $\Theta_f$, $\Omega_f\gg (b-a)$  and $\Omega_g$, for which the derivatives satisfy
\begin{align}
\label{cond-1}
f^{(i)}(x)\ll \Theta_f/\Omega_f^i,\;\;\;\;\;\;\;g^{(j)}(x)\ll 1/\Omega_g^j.
\end{align}
For the second assertion we will moreover require
\begin{align}
\label{cond-2}
f^{(2)}(x)\gg \Theta_f/\Omega_f^2.
\end{align}

\begin{lemma}
\label{sp}
Suppose $f$ and $g$ are smooth real valued satisfying \eqref{cond-1} for $i=2,3$ and $j=0,1,2$. Suppose $g(a)=g(b)=0$.

(1) Suppose $f'$ and $f''$ do not vanish in $[a,b]$. Let $\Lambda=\min_{[a,b]}|f'(x)|$. Then we have 
$$
\mathfrak{I}\ll \frac{\Theta_f}{\Omega_f^2\Lambda^3}\left(1+\frac{\Omega_f}{\Omega_g}+\frac{\Omega_f^2}{\Omega_g^2}\frac{\Lambda}{\Theta_f/\Omega_f}\right).
$$

(2) Suppose $f'$ changes sign from negative to positive at the unique point $x_0\in (a,b)$. Let $\kappa=\min\{b-x_0,x_0-a\}$. Further suppose \eqref{cond-1} holds for $i=4$ and \eqref{cond-2} holds. Then we have 
$$
\mathfrak{I}= \frac{g(x_0)e(f(x_0)+1/8)}{\sqrt{f''(x_0)}}+O\left(\frac{\Omega_f^4}{\Theta_f^2\kappa^3}+\frac{\Omega_f}{\Theta_f^{3/2}}+\frac{\Omega_f^3}{\Theta_f^{3/2}\Omega_g^2}\right).
$$\\
\end{lemma} 

Finally we recall the second derivative bound for exponential integrals in two variables. Let
\begin{align}
\label{double-expo}
\mathfrak{I}_{(2)}=\int_{a}^b\int_c^d g(x,y)e(f(x,y))\mathrm{d}y\mathrm{d}x
\end{align}
where $f$ and $g$ are smooth real valued functions. First suppose $g=1$, and we have positive parameters $r_1$ and $r_2$ such that in the rectangle $[a,b]\times [c,d]$ we have
\begin{align}
\label{exp-int-2-var}
f^{(2,0)}(x,y)\gg r_1^2,\;\;\;\;f^{(0,2)}(x,y)\gg r_2^2,\;\;\;\;f^{(2,0)}(x,y)f^{(0,2)}(x,y)-\left[f^{(1,1)}(x,y)\right]^2\gg r_1^2r_2^2,
\end{align}
where $f^{(i,j)}(x,y)=\frac{\partial^{i+j}}{\partial x^i\partial y^j}f(x,y)$ and the implied constants are absolute. Then we have (see \cite{S})
$$
\mathfrak{I}_{(2)}\ll \frac{1}{r_1r_2}.
$$ 
To extend this result to smooth $g$ with $\mathrm{Supp}(g)\subset (a,b)\times (c,d)$, we apply integration by parts once in each variable. To state the result we define the total variation of $g$ to be 
$$
\text{var}(g):=\int_{a}^b\int_c^d\left|\frac{\partial^2}{\partial  x\partial y}g(x,y)\right|\mathrm{d}y\mathrm{d}x.
$$ 
\begin{lemma}
\label{exp-2-var-lemma}
Suppose $f$, $g$, $r_1$ and $r_2$ are as above satisfying the condition \eqref{exp-int-2-var}. Then we have
$$
\mathfrak{I}_{(2)}\ll \frac{\mathrm{var}(g)}{r_1r_2}
$$ 
with an absolute implied constant.\\
\end{lemma}

\subsection{An integral}
\label{anintegral}

Let $W$ be a smooth real valued function with $\mathrm{Supp}(W)\subset [a,b]\subset (0,\infty)$. Suppose $W^{(j)}(x)\ll_{a,b,j}1$. Consider the exponential integral
\begin{align}
\label{int-tran}
W^\dagger(r,s)=\int_0^\infty W(x)e(-rx)x^{s-1}\mathrm{d}x
\end{align}
where $r\in \mathbb R$ and $s=\sigma+i\beta \in \mathbb C$. In particular $W^\dagger(r,1)$ is the Fourier transform of $W$ and $W^\dagger(0,s)$ is the Mellin transform of $W$. The integral is of the form \eqref{exp-integral} with 
$$
g(x)=W(x)x^{\sigma-1}\;\;\;\;\text{and}\;\;\;\;f(x)=-rx+\frac{1}{2\pi}\beta \log x.
$$
Then 
$$
f'(x)=-r+\frac{1}{2\pi}\frac{\beta}{x}\;\;\;\;\text{and}\;\;\;\;f^{(j)}(x)=(-1)^{j-1}(j-1)!\frac{1}{2\pi}\frac{\beta}{x^j}
$$
for $j\geq 2$. The unique stationary point is given by 
$$
x_0=\frac{\beta}{2\pi r},
$$
and we can write
$$
f'(x)=\frac{\beta}{2\pi}\left(\frac{1}{x}-\frac{1}{x_0}\right)=r\left(\frac{x_0}{x}-1\right).
$$
Suppose $a$, $b$ and $\sigma$ are fixed and we are interested in the dependence of the integral on $\beta$ and $r$. If $x_0\notin [a/2,2b]$ then in the support of the integral, i.e. $W(x)\neq 0$, we have $|f'(x)|\gg \max\{|\beta|,|r|\}$ and $|f^{(j)}(x)|\ll_j |\beta|$, where the implied constants depend on $a$, $b$ and $\sigma$. So in this case $W^\dagger(r,s)\ll_j \min\{|\beta|^{-j},|r|^{-j}\}$, where again the implied constant depends on $a$, $b$ and $\sigma$. On the other hand if $x_0\in [a/2,2b]$ then using the second statement of Lemma~\ref{sp} (with $\Theta_f=|\beta|$ and $\Omega_f=\Omega_g=1$) we get
$$
W^\dagger(r,s)= \frac{g(x_0)e(f(x_0)+1/8)}{\sqrt{f''(x_0)}}+O\left(|\beta|^{-3/2}\right).
$$
The error term can also be written as $O(|r|^{-3/2})$, as $x_0\in [a/2,2b]$ implies that $|r|\asymp |\beta|$. Note that for $\beta>0$ we need to take conjugate so that the conditions of the lemma are satisfied. Also we note that the above asymptotic holds regardless the location of $x_0$. For the following statement we take $\sqrt{-1}=e^{\pi i/2}$.
\begin{lemma}
\label{the-integral}
Let $W$ be a smooth real valued function with $\mathrm{Supp}(W)\subset [a,b]\subset (0,\infty)$ and $W^{(j)}(x)\ll_{a,b,j}1$. Let $r\in \mathbb R$ and $s=\sigma+i\beta \in \mathbb C$. We have
$$
W^\dagger(r,s)= \frac{\sqrt{2\pi}e(1/8)}{\sqrt{-\beta}}W\left(\frac{\beta}{2\pi r}\right)\left(\frac{\beta}{2\pi r}\right)^\sigma\left(\frac{\beta}{2\pi er}\right)^{i\beta}+O\left(\min\{|\beta|^{-3/2},|r|^{-3/2}\}\right),
$$
where the implied constant depends on $a$, $b$ and $\sigma$. We also have
$$
W^\dagger(r,s)=O_{a,b,\sigma,j}\left(\min\left\{\left(\frac{1+|\beta|}{|r|}\right)^{j}, \left(\frac{1+
|r|}{|\beta|}\right)^{j}\right\}\right).
$$ 
\end{lemma}



\section{Application of summation formula}
\label{pvs}

\subsection{Applying Poisson summation}\label{poi-1}
For simplicity let us assume that $t>2$. First we will apply the Poisson summation formula on the sum over $m$ in \eqref{S(N)-cm}, i.e.
$$
\mathop{\sum}_{m=1}^\infty m^{-i(t+v)}e\left(-\frac{m\bar a}{q}\right)e\left(\frac{mx}{aq}\right)U\left(\frac{m}{N}\right).
$$
Breaking the sum into congruence classes modulo $q$ we get
$$
\sum_{\alpha\bmod{q}}e\left(-\frac{\alpha\bar a}{q}\right)\mathop{\sum}_{m\in\mathbb Z} (\alpha+mq)^{-i(t+v)}e\left(\frac{(\alpha+mq)x}{aq}\right)U\left(\frac{\alpha+mq}{N}\right).
$$
Then applying Poisson we obtain
$$
\sum_{\alpha\bmod{q}}e\left(-\frac{\alpha\bar a}{q}\right)\mathop{\sum}_{m\in\mathbb Z} \int_{\mathbb R}(\alpha+yq)^{-i(t+v)}e\left(\frac{(\alpha+yq)x}{aq}\right)U\left(\frac{\alpha+yq}{N}\right)e(-my)\mathrm{d}y.
$$
Making the change of variable $(\alpha+yq)/N\mapsto u$ and executing the resulting complete character sum we arrive at
\begin{align}
\label{a-det}
N^{1-i(t+v)}\mathop{\sum}_{\substack{m\in\mathbb Z\\m\equiv \bar a\bmod{q}}} \int_{\mathbb R}U\left(u\right)u^{-i(t+v)}e\left(\frac{N(x-ma)}{aq}u\right)\mathrm{d}u.
\end{align}
The above integral, in the notation of Subsection~\ref{anintegral} is 
\begin{align}
\label{I1}
U^{ \dagger}\left(N(aq)^{-1}(ma-x), 1-i(t+v)\right).
\end{align}
Recall that $a\asymp (N/K)^{1/2}$, and by our choice, see \eqref{cond-K}, $K<t^{2-\varepsilon}/N$. So $|N(aq)^{-1}(ma-x)|\asymp q^{-1}N|m|$ if $m\neq 0$, and $|N(aq)^{-1}(ma-x)|\ll q^{-1}(NK)^{1/2}$ if $m=0$. Applying the second statement of Lemma~\ref{the-integral} it follows that the contribution of the zero frequency $m=0$ (which occurs only for $q=1$ due to the condition $(m,q)=1$) in \eqref{a-det} is negligibly small, and also the contribution of the tail $|m|\gg qt^{1+\varepsilon}/N$ is negligibly small. We only need to consider $m$ with $1\leq |m|\ll qt^{1+\varepsilon}/N$, which in turn implies that we only need to focus on $q$ which are in the range 
$$
N/t^{1+\varepsilon}\ll q \leq Q=(N/K)^{1/2}.
$$ 
Taking a dyadic subdivision we conclude the following. 

\begin{lemma}
\label{pre-lemma}
Suppose $N$ and $K$ satisfy \eqref{cond-K}, then we have
\begin{align}
\label{dyadic}
S^+(N)=\frac{N}{K}\sum_{\substack{N/t^{1+\varepsilon}\ll C \leq (N/K)^{1/2}\\ C\;\mathrm{dyadic}}} S(N,C)+O(t^{-2012})
\end{align}
where
\begin{align}
\label{SNC}
S(N,C)=\int_0^1&\int_{\mathbb R}N^{-i(t+v)}V\left(\frac{v}{K}\right)\mathop{\sum}_{C< q\leq 2C}\mathop{\sum}_{\substack{(m,q)=1\\1\leq |m|\ll qt^{1+\varepsilon}/N}}\frac{1}{aq}\:U^{ \dagger}\left(\frac{N(ma-x)}{aq}, 1-i(t+v)\right)\\
\nonumber &\times \mathop{\sum}_{n=1}^\infty \lambda(1,n)e\left(\frac{nm}{q}\right)n^{iv}e\left(-\frac{nx}{aq}\right)V\left(\frac{n}{N}\right)\mathrm{d}v\mathrm{d}x.
\end{align}
Here $a=a_Q(m,q)$ is the unique multiplicative inverse of $m$ modulo $q$ in the range $(Q,q+Q]$.
\end{lemma}


\subsection{Applying Voronoi summation}
Applying Lemma~\ref{lem-voronoi} we get
\begin{align}
\label{vor}
\mathop{\sum}_{n=1}^\infty & \lambda(1,n)e\left(\frac{nm}{q}\right)n^{iv}e\left(-\frac{nx}{aq}\right)V\left(\frac{n}{N}\right)\\
\nonumber =&qN^{iv}\sum_\pm \sum_{n_1|q}\sum_{n_2=1}^\infty \frac{\lambda(n_2,n_1)}{n_1n_2}S(\bar m, \pm n_2; q/n_1)\mathfrak{I}\left(\frac{n_1^2n_2}{q^3},\frac{x}{aq}\right)
\end{align}
where
\begin{align*}
\mathfrak{I}(r,r')=\frac{1}{2\pi i}\int_{(\sigma)}\left(rN\right)^{-s}\gamma_\pm(s)V^\dagger(Nr',-s+iv)\mathrm{d}s.
\end{align*}
(Here $V^\dagger$ is as defined in Subsection~\ref{anintegral}.) Using Stirling approximation we get that
\begin{align}
\label{stir}
\left|\gamma_\pm (s)\right|\ll_{\pi,\sigma} 1+|\tau|^{3\sigma+3/2}
\end{align}
where $s=\sigma+i\tau$ and $\sigma\geq -1/2$. Also from the second statement of Lemma~\ref{the-integral} we get that 
$$
V^\dagger\left(\frac{Nx}{aq}, -s+iv\right)\ll_{\sigma,j} \min\left\{1,\left(\frac{NK^{1/2}}{q|\tau-v|}\right)^j\right\}.
$$
Shifting the contour to $\sigma=M$ (a large positive integer) and taking $j=3M+3$ we can make the integral in \eqref{vor} negligibly small if
$
n_1^2n_2\gg N^{1/2}K^{3/2}t^\varepsilon.
$
For smaller values of $n_1^2n_2$ we move the contour to $\sigma=-1/2$, and obtain
\begin{align*}
\mathfrak{I}\left(\frac{n_1^2n_2}{q^3},\frac{x}{aq}\right)=\frac{1}{2\pi}\left(\frac{n_1^2n_2N}{q^3}\right)^{1/2}&\sum_{J\in\mathcal J}\mathop{\int}_{\mathbb R}\left(\frac{n_1^2n_2N}{q^3}\right)^{-i\tau}\gamma_\pm\left(-\frac{1}{2}+i\tau\right)\\
&\times V^\dagger\left(\frac{Nx}{aq}, \frac{1}{2}-i\tau+iv\right)W_J(\tau)\mathrm{d}\tau
+O(t^{-20120}).
\end{align*} 
Here $\mathcal J$ is a collection of $O(\log t)$ many real numbers in the interval $[-(NK)^{1/2}t^\varepsilon/C, (NK)^{1/2}t^\varepsilon/C]$, containing $0$. For each $J$ we have a smooth bump function (non-negative real valued) $W_J(x)$ satisfying $x^\ell W_J^{(\ell)}(x)\ll_\ell 1$ for all $\ell\geq 0$. For $J=0$ the function $W_0(x)$ is supported in $[-1,1]$ and satisfies the stricter bound $W_J^{(\ell)}(x)\ll_\ell 1$. For each $J> 0$ (resp. $J<0$) the function $W_J(x)$ is supported in $[J,4J/3]$ (resp. $[4J/3,J]$).  Finally we require that
$$
\sum_{J\in\mathcal J}W_J(x)=1,\;\;\;\;\text{for}\;\;x\in [-(NK)^{1/2}t^\varepsilon/C, (NK)^{1/2}t^\varepsilon/C].
$$
In short the collection $W_J$ is a smooth partition of unity. The precise definition of the function or the collection will not be required. \\

\begin{lemma}
\label{lemma-snc}
Let $N$ and $K$ satisfy \eqref{cond-K}, and suppose $N/t^{1+\varepsilon}\ll C\ll (N/K)^{1/2}$. We have
\begin{align*}
S(N,C)=\frac{N^{1/2-it}K}{2\pi}&\sum_\pm \sum_{J\in\mathcal J} \mathop{\sum\sum}_{n_1^2n_2\ll N^{1/2}K^{3/2}t^\varepsilon} \frac{\lambda(n_2,n_1)}{n_2^{1/2}}\\
&\times \mathop{\sum\sum}_{\substack{C< q\leq 2C,\; (m,q)=1\\1\leq |m|\ll qt^{1+\varepsilon}/N\\n_1|q}}\frac{S(\bar m, \pm n_2; q/n_1)}{aq^{3/2}}\mathfrak{I}_{\pm}^\star(q,m,n_1^2n_2)+O(t^{-2012}),
\end{align*}
where
\begin{align*}
\mathfrak{I}_{\pm}^\star(q,m,n)=\int_{\mathbb R}\left(\frac{nN}{q^3}\right)^{-i\tau}\gamma_\pm\left(-\frac{1}{2}+i\tau\right)\mathfrak{I}^{\star\star}(q,m,\tau)W_J(\tau)\mathrm{d}\tau,
\end{align*}
and
\begin{align*}
\mathfrak{I}^{\star\star}(q,m,\tau)=\int_0^1\mathop{\int}_{\mathbb R}V\left(v\right)V^\dagger\left(\frac{Nx}{aq}, \frac{1}{2}-i\tau+iKv\right)U^{ \dagger}\left(\frac{N(ma-x)}{aq}, 1-i(t+Kv)\right)\mathrm{d}v\mathrm{d}x.
\end{align*}
\end{lemma}

In the next section we will analyse the integrals further.


\section{Analysis of the integrals}
\label{aoi}

The aim of this section is to prove Lemma~\ref{i3-final} 
which gives a decomposition of $\mathfrak{I}^{\star\star}(q,m,\tau)$ for $|\tau|\ll K^{1/2}t^{1+\varepsilon}/N^{1/2}$.\\

\subsection{Stationary phase analysis for $V^\dagger$ and $U^{\dagger}$} We apply the first statement of Lemma~\ref{the-integral} 
to conclude that
\begin{align*}
U^{\dagger}\left(\frac{N(ma-x)}{aq}, 1-i(t+Kv)\right)=e^{\pi i/4}&\frac{(t+Kv)^{1/2}\:aq}{(2\pi)^{1/2}N(x-ma)}U\left(\frac{(t+Kv)aq}{2\pi N(x-ma)}\right)\\
\times &\left(\frac{(t+Kv)aq}{2\pi eN(x-ma)}\right)^{-i(t+Kv)}+O\left(\frac{1}{t^{3/2}}\right).
\end{align*}
The error term makes a contribution of size $O\left(t^{\varepsilon-3/2}\right)$ towards $\mathfrak{I}^{\star\star}(q,m,\tau)$, and we get 
\begin{align}
\label{after-u-int}
\mathfrak{I}^{\star\star}(q,m,\tau)=c_1\frac{aq}{N}\int_0^1&\int_{\mathbb R}V(v)V^\dagger\left(\frac{Nx}{aq}, \frac{1}{2}-i\tau+iKv\right)\frac{(t+Kv)^{1/2}}{(x-ma)}U\left(\frac{(t+Kv)aq}{2\pi N(x-ma)}\right)\\
\nonumber &\times \left(\frac{(t+Kv)aq}{2\pi eN(x-ma)}\right)^{-i(t+Kv)}\mathrm{d}v\mathrm{d}x+O\left(\frac{t^{\varepsilon}}{t^{3/2}}\right)
\end{align}
for some constant $c_1$, and an absolute implied constant.\\

Next we study the integral 
$
V^\dagger\left(Nx/aq, 1/2-i\tau+iKv\right)
$
using Lemma~\ref{the-integral}. 
We have
\begin{align*}
V^\dagger\left(\frac{Nx}{aq}, \frac{1}{2}-i\tau+iKv\right)=e^{-\pi i/4}\frac{(aq)^{1/2}}{(2Nx)^{1/2}}&V\left(\frac{(Kv-\tau)aq}{2\pi Nx}\right)\left(\frac{(Kv-\tau)aq}{2\pi eNx}\right)^{i(Kv-\tau)}\\
&+O\left(\min\left\{\left(\frac{aq}{Nx}\right)^{3/2},\frac{1}{|\tau-Kv|^{3/2}}\right\}\right).
\end{align*}
Hence up to a constant $\mathfrak{I}^{\star\star}(q,m,\tau)$, is given by
\begin{align*}
\left(\frac{aq}{N}\right)^{3/2}&\int_0^1\int_\mathbb{R} V\left(v\right)\frac{(t+Kv)^{1/2}}{x^{1/2}(x-ma)}U\left(\frac{(t+Kv)aq}{2\pi N(x-ma)}\right)V\left(\frac{(Kv-\tau)aq}{2\pi Nx}\right)\\
&\times \left(\frac{(t+Kv)aq}{2\pi eN(x-ma)}\right)^{-i(t+Kv)}\left(\frac{(Kv-\tau)aq}{2\pi eNx}\right)^{i(Kv-\tau)}\mathrm{d}v\mathrm{d}x+O(E^{\star\star}+t^{\varepsilon-3/2})
\end{align*}
where the error term is given by (since $uU(u)\ll 1$)
\begin{align*}
E^{\star\star}=\frac{1}{t^{1/2}}\int_0^1&\int_{1}^2\:\min\left\{\left(\frac{aq}{Nx}\right)^{3/2},\frac{1}{|\tau-Kv|^{3/2}}\right\}\mathrm{d}v\mathrm{d}x.
\end{align*}
To estimate the error term we split the inner integral over $v$ into pieces. Indeed the first term in the integrand is smaller than the second term if and only if $v$ lies in the interval
$$
\frac{\tau}{K}-\frac{Nx}{aqK}<v<\frac{\tau}{K}+\frac{Nx}{aqK}.
$$ 
If $|\tau|\geq 10K$, this interval does not intersect $[1,2]$ unless $Nx/aq\asymp |\tau|$. In this case we use the trivial bound $O(1)$ for the inner integral over $v$. On the other hand if $|\tau|<10K$ the inner integral is bounded by $Nx/aqK$. Hence the contribution of the case where the first term in the integrand is smaller in $E^{\star\star}$ is bounded by
$$
\frac{1}{t^{1/2}}\int_0^1 \left(\frac{aq}{Nx}\right)^{1/2}\:\frac{1}{K}\:\mathbf{1}_{|\tau|<10K}\mathrm{d}x+\frac{1}{t^{1/2}}\int_0^1 \left(\frac{aq}{Nx}\right)^{1/2}\:\frac{1}{|\tau|}\:\mathbf{1}_{|\tau|\geq 10K}\mathrm{d}x.
$$
(Here $\mathbf{1}_s$ is the indicator function of the statement $s$, i.e. it takes the value $1$ if the statement $s$ is true, and vanishes otherwise.) This is dominated by
\begin{align}
\label{bd-revised}
O\left(\frac{1}{t^{1/2}K^{3/2}}\min\left\{1,\frac{10K}{|\tau|}\right\}t^{\varepsilon}\right).
\end{align}\\

Next we estimate the contribution of the case where the second term in the integrand in $E^{\star\star}$ is smaller. This is given by
\begin{align*}
\frac{1}{t^{1/2}}\mathop{\int_0^1\int_{1}^2}_{|\tau-Kv|>Nx/aq}\:\frac{1}{|\tau-Kv|^{3/2}}\mathrm{d}v\mathrm{d}x &\ll \frac{1}{t^{1/2}}\int_0^1 \left(\frac{aq}{Nx}\right)^{1/2+\varepsilon}\:\int_0^1\:\frac{1}{|\tau-Kv|^{1-\varepsilon}}\mathrm{d}v\mathrm{d}x\\
&\ll t^{\varepsilon}\:\frac{Q}{t^{1/2}N^{1/2}}\:\frac{1}{K}\min\left\{1,\frac{10K}{|\tau|}\right\}.
\end{align*}
This is again dominated by the bound given in \eqref{bd-revised}. 
Also it follows from (5), and the inequality $|\tau|\ll K^{1/2}t^{1+\varepsilon}/N^{1/2}$, that $t\gg K^{3/2}\max\{1,|\tau|/10K\}$.
Hence
$$
E^{\star\star}+t^{-3/2+\varepsilon}\ll \frac{1}{t^{1/2}K^{3/2}}\min\left\{1,\frac{10K}{|\tau|}\right\}t^{\varepsilon},
$$
and we conclude that
\begin{align}
\label{i3-mid}
\mathfrak{I}^{\star\star}(q,m,\tau)=c_2\left(\frac{aq}{N}\right)^{1/2}&\int_0^1\int_{\mathbb R}V\left(v\right)\frac{(t+Kv)^{1/2}aq}{x^{1/2}(x-ma)N}U\left(\frac{(t+Kv)aq}{2\pi N(x-ma)}\right)V\left(\frac{(Kv-\tau)aq}{2\pi Nx}\right)\\
\nonumber &\times \left(\frac{(t+Kv)aq}{2\pi eN(x-ma)}\right)^{-i(t+Kv)}\left(\frac{(Kv-\tau)aq}{2\pi eNx}\right)^{i(Kv-\tau)}\mathrm{d}v\mathrm{d}x\\
\nonumber &+O\left(\frac{1}{t^{1/2}K^{3/2}}\min\left\{1,\frac{10K}{|\tau|}\right\}t^{\varepsilon}\right)
\end{align}
for some absolute constant $c_2$.\\

\subsection{Integral over $v$} Now we will study the integral over $v$ in \eqref{i3-mid}. This term vanishes unless $m<0$. For $x<1/K$ we bound the integral trivially. Indeed, in this case the weight function restricts the integral over $v$ to a range of length $N/K^2aq$. (Recall that $u^rU(u)\ll_r 1$ and $u^rV(u)\ll_r 1$ for any $r\in \mathbb{R}$.) So estimating trivially we get
\begin{align}
\label{small-x}
\left(\frac{aq}{N}\right)^{1/2}\int_0^{1/K}\int_{\mathbb R}V\left(v\right)&\frac{(t+Kv)^{1/2}\:aq}{x^{1/2}(x-ma)N}\:U\left(\frac{(t+Kv)aq}{2\pi N(x-ma)}\right)V\left(\frac{(Kv-\tau)aq}{2\pi Nx}\right)\mathrm{d}v\mathrm{d}x\\
\nonumber &\ll \frac{1}{t^{1/2}\:K^{5/2}}\left(\frac{N}{aq}\right)^{1/2}.
\end{align}
Let us now take $x\in [1/K,1]$. Temporarily we set
$$
f(v)=-\frac{t+Kv}{2\pi}\log \left(\frac{(t+Kv)aq}{2\pi eN(x-ma)}\right)+\frac{Kv-\tau}{2\pi}\log \left(\frac{(Kv-\tau)aq}{2\pi eNx}\right)
$$
and
$$
g(v)=\frac{t^{1/2}(t+Kv)^{1/2}\:aq}{N(x-ma)}V\left(v\right)U\left(\frac{(t+Kv)aq}{2\pi N(x-ma)}\right)V\left(\frac{(Kv-\tau)aq}{2\pi Nx}\right).
$$
We are multiplying by an extra $t^{1/2}$ to balance the size of the function. Then
\begin{align*}
f'(v)=-\frac{K}{2\pi}\log \left(\frac{(t+Kv)x}{(Kv-\tau)(x-ma)}\right),\;\;\;\;\;\;f^{(j)}(v)=-\frac{(j-1)!(-K)^j}{2\pi(t+Kv)^{j-1}}+\frac{(j-1)!(-K)^j}{2\pi (Kv-\tau)^{j-1}},
\end{align*}
for $j\geq 2$. The stationary phase is given by
\begin{align*}
v_0=-\frac{(t+\tau)x-\tau ma}{Kma}.
\end{align*}
In the support of the integral we have
$$
f^{(j)}\asymp \frac{Nx}{aq}\left(\frac{Kaq}{Nx}\right)^j
$$
for $j\geq 2$, and 
\begin{align*}
g^{(j)}(v)\ll \left(1+\frac{Kaq}{Nx}\right)^j.
\end{align*}
for $j\geq 0$. Moreover we can write
$$
f'(v)=\frac{K}{2\pi}\log \left(1+\frac{K(v_0-v)}{(t+Kv)}\right)-\frac{K}{2\pi}\log \left(1+\frac{K(v_0-v)}{(Kv-\tau)}\right).
$$
In the support of the integral we have $0\leq Kv-\tau\ll N/aq \ll K^{1/2}t^{1+\varepsilon}/N^{1/2}\ll t^{3/4+\varepsilon}$ (recall that $t<N$, $K<t^{1/2}$ and $N/t^{1+\varepsilon}<q$). It follows that if $v_0\notin [.5,3]$ then in the support of the integral we have
$$
|f'(v)|\gg K^{1-\varepsilon}\min\left\{1,\frac{Kaq}{Nx}\right\}.
$$
Applying the first statement of Lemma~\ref{sp} with
$$
\Theta_f=\frac{Nx}{aq},\;\;\;\;\Omega_f=\frac{Nx}{Kaq},\;\;\;\;\Omega_g=\min\left\{1,\frac{Nx}{Kaq}\right\},\;\;\;\text{and}\;\;\;\Lambda=K^{1-\varepsilon}\min\left\{1,\frac{Kaq}{Nx}\right\}
$$
we obtain the bound
\begin{align}
\label{v-int-bd}
\int_{\mathbb R}g(v)e(f(v))\mathrm{d}v\ll  \frac{\Theta_f}{\Omega_f^2\Lambda^3}\left(1+\frac{\Omega_f}{\Omega_g}+\frac{\Omega_f^2}{\Omega_g^2}\frac{\Lambda}{\Theta_f/\Omega_f}\right)t^\varepsilon.
\end{align}
On the other hand if $v_0\in [.5,3]$ then treating the integral as a finite integral over the range $[.1,4]$ and applying the second part of Lemma~\ref{sp} it follows that
\begin{align}
\label{v-int-bd2}
\int_{\mathbb R}g(v)e(f(v))\mathrm{d}v= \frac{g(v_0)e(f(v_0)+1/8)}{\sqrt{f''(v_0)}}+O\left(\left(\frac{\Omega_f^4}{\Theta_f^2}+\frac{\Omega_f}{\Theta_f^{3/2}}+\frac{\Omega_f^3}{\Theta_f^{3/2}\Omega_g^2}\right)t^\varepsilon\right).
\end{align}
Notice that we have $\kappa>.4$. The bound from \eqref{v-int-bd} and the error term of \eqref{v-int-bd2} together make a total contribution of size
\begin{align}
\label{omega-theta-bd}
O\left(\frac{1}{t^{1/2}}\frac{N}{aqK^3}t^\varepsilon\right)
\end{align}
in \eqref{i3-mid}. We arrive at this through explicit calculation (estimating the integral over $x$ trivially).  If $x\leq Kaq/N$ then 
$
\Omega_g=\Omega_f,
$
and $\Lambda=K^{1-\varepsilon}$. We find that the error term in \eqref{v-int-bd2} is dominated by $O((aq)^{1/2}t^\varepsilon/K(Nx)^{1/2})$ and the bound in \eqref{v-int-bd} is dominated by $O(aqt^\varepsilon/NKx)$. Observe that as $x>1/K$ we have $aq/NKx\ll (aq)^{1/2}/K(Nx)^{1/2}$. Hence the total contribution of this part is dominated by 
\begin{align*}
t^\varepsilon\left(\frac{aq}{tN}\right)^{1/2}\int_{1/K}^{\max\{1/K,Kaq/N\}} \frac{(aq)^{1/2}}{K(Nx)^{1/2}}\mathrm{d}x\ll t^\varepsilon\frac{(aq)^{3/2}}{t^{1/2}N^{3/2}K^{1/2}}.
\end{align*}
This is dominated by \eqref{omega-theta-bd} as $aq\ll N/K$. On the other hand if $x> Kaq/N$ then $\Omega_g=1$, and $\Lambda=K^{1-\varepsilon}/\Omega_f$. So the error term in \eqref{v-int-bd2} is dominated by $O(N^{3/2}t^\varepsilon/(K^2aq)^{3/2})$ and the bound in \eqref{v-int-bd} is dominated by $O(N^3t^\varepsilon/K^5(aq)^3)$. Now the lower bound on $K$ as given in \eqref{cond-K} implies that $N^3/K^5(aq)^3\ll N^{3/2}t^\varepsilon/(K^2aq)^{3/2}$ for $q>N/t^{1+\varepsilon}$. Hence the total contribution of this part is dominated by 
\begin{align*}
t^\varepsilon\left(\frac{aq}{tN}\right)^{1/2}\int_0^1 \frac{N^{3/2}}{(K^2aq)^{3/2}}\mathrm{d}x\ll t^\varepsilon\frac{N}{t^{1/2}K^3aq}.
\end{align*}
Also the above term dominates the bound obtained in \eqref{small-x}.\\ 

We conclude that
\begin{align*}
\left(\frac{aq}{tN}\right)^{1/2}\int_0^1 \frac{1}{x^{1/2}}\:\int_\mathbb{R} g(v)e(f(v))\mathrm{d}v\:\mathrm{d}x=&\left(\frac{aq}{tN}\right)^{1/2}\int_{1/K}^1 \frac{1}{x^{1/2}}\frac{g(v_0)e(f(v_0)+1/8)}{\sqrt{f''(v_0)}}\mathrm{d}x\\
&+O\left(\frac{1}{t^{1/2}}\frac{N}{aqK^3}t^\varepsilon\right).
\end{align*}
By explicit computation we get
$$
f(v_0)=-\frac{t+\tau}{2\pi}\log\left(\frac{-(t+\tau)q}{2\pi eNm}\right),\;\;\;\;f''(v_0)=\frac{K^2(ma)^2}{2\pi (t+\tau)(x-ma)x}
$$
and
$$
g(v_0)=\frac{aq}{N}\left(\frac{-t(t+\tau)}{ma(x-ma)}\right)^{1/2}V\left(\frac{\tau}{K}-\frac{(t+\tau)x}{Kma}\right) U\left(\frac{-(t+\tau)q}{2\pi Nm}\right)V\left(\frac{-(t+\tau)q}{2\pi Nm}\right).
$$
Plugging these in, and using the fact that $U(y)V(y)=V(y)$, the leading term in the above expression reduces to
\begin{align*}
c_3\frac{t+\tau}{K}\left(\frac{q}{-mN}\right)^{3/2}\:V\left(-\frac{(t+\tau)q}{2\pi Nm}\right)\left(-\frac{(t+\tau)q}{2\pi eNm}\right)^{-i(t+\tau)}\:\int_{1/K}^1V\left(\frac{\tau}{K}-\frac{(t+\tau)x}{Kma}\right)\mathrm{d}x
\end{align*}
for some absolute constant $c_3$. We can now extend the integral to the range $[0,1]$ at a cost of an error term which is dominated by the error term in \eqref{i3-mid}.\\

Now we will summarize the main output of this section. Set
\begin{align}
\label{bctau}
B(C,\tau)=\frac{1}{t^{1/2}K^{3/2}}\min\left\{1,\frac{10K}{|\tau|}\right\}+\frac{1}{t^{1/2}K^{5/2}}\frac{N^{1/2}}{C}.
\end{align}
Note that
\begin{align}
\label{bctau-int}
\int_{-\frac{(NK)^{1/2}}{C}t^\varepsilon}^{+\frac{(NK)^{1/2}}{C}t^\varepsilon}B(C,\tau)\mathrm{d}\tau\ll \frac{t^\varepsilon}{t^{1/2}K^{1/2}}\left\{1+\frac{N}{C^2K^{3/2}}\right\}.
\end{align}\\

\begin{lemma}
\label{i3-final}
Suppose $C\leq q\leq 2C$, with $N/t^{1+\varepsilon}\ll C \leq (N/K)^{1/2}$, and $N$, $K$ satisfy \eqref{cond-K}. Suppose $t>2$ and $|\tau|\ll K^{1/2}t^{1+\varepsilon}/N^{1/2}$. We have 
\begin{align*}
\mathfrak{I}^{\star\star}(q,m,\tau)=\mathfrak{J}_{1}(q,m,\tau)+\mathfrak{J}_{2}(q,m,\tau)
\end{align*}
where
\begin{align*}
\mathfrak{J}_{1}(q,m,\tau)=\frac{c_4}{(t+\tau)^{1/2}\:K}\left(-\frac{(t+\tau)q}{2\pi eNm}\right)^{3/2-i(t+\tau)}V\left(-\frac{(t+\tau)q}{2\pi Nm}\right)\int_0^1V\left(\frac{\tau}{K}-\frac{(t+\tau)x}{Kma}\right)\mathrm{d}x
\end{align*}
for some absolute constant $c_4$ and
\begin{align*}
\mathfrak{J}_{2}(q,m,\tau):=\mathfrak{I}^{\star\star}(q,m,\tau)-\mathfrak{J}_{1}(q,m,\tau)=O\left(B(C,\tau)t^\varepsilon\right),
\end{align*}
with $B(C,\tau)$ as defined in \eqref{bctau}.
\end{lemma}

Consequently from Lemma~\ref{lemma-snc} we derive the following decomposition for $S(N,C)$.

\begin{lemma}
\label{pre-cauchy}
We have
$$
S(N,C)=\sum_{J\in\mathcal J}\left\{S_{1,J}(N,C)+S_{2,J}(N,C)\right\}+O(t^{-2012})
$$
where
\begin{align*}
S_{\ell,J}(N,C)=\frac{N^{1/2-it}K}{2\pi}&\sum_\pm \mathop{\sum\sum}_{n_1^2n_2\ll N^{1/2}K^{3/2}t^\varepsilon} \frac{\lambda(n_2,n_1)}{n_2^{1/2}}\\
&\times \mathop{\sum\sum}_{\substack{C< q\leq 2C,\; (m,q)=1\\1\leq |m|\ll qt^{1+\varepsilon}/N\\n_1|q}}\frac{S(\bar m, \pm n_2; q/n_1)}{aq^{3/2}}\mathfrak{I}_{\ell,J,\pm}(q,m,n_1^2n_2),
\end{align*}
and
\begin{align*}
\mathfrak{I}_{\ell,J,\pm}(q,m,n)=\int_{\mathbb R}\left(\frac{nN}{q^3}\right)^{-i\tau}\gamma_\pm\left(-\frac{1}{2}+i\tau\right)\mathfrak{J}_{\ell}(q,m,\tau)W_J(\tau)\mathrm{d}\tau
\end{align*}
with $\mathfrak{J}_{\ell}(q,m,\tau)$ as defined in Lemma~\ref{i3-final}.
\end{lemma}


\section{Application of Cauchy inequality and Poisson summation - I}

In this section we will estimate 
$$
S_2(N,C):=\sum_{J\in\mathcal J}S_{2,J}(N,C).
$$ 
In this case we will not need any cancellation in the integral over $\tau$. 

\subsection{Applying Cauchy inequality and Poisson summation} Taking a dyadic segmentation, and using the bound $|\gamma_\pm\left(-\frac{1}{2}+i\tau\right)|\ll 1$, we get
\begin{align*}
S_2(N,C)\leq t^\varepsilon N^{1/2}K\int_{-\frac{(NK)^{1/2}}{C}t^\varepsilon}^{+\frac{(NK)^{1/2}}{C}t^\varepsilon}\sum_\pm &\sum_{\substack{1\leq L\ll N^{1/2}K^{3/2}t^\varepsilon\\ \text{dyadic}}} \mathop{\sum\sum}_{n_1, n_2} \frac{|\lambda(n_2,n_1)|}{n_2^{1/2}}U\left(\frac{n_1^2n_2}{L}\right)\\
&\times \Bigl|\mathop{\sum\sum}_{\substack{C< q\leq 2C,\; (m,q)=1\\1\leq |m|\ll qt^{1+\varepsilon}/N\\n_1|q}}\frac{S(\bar m, \pm n_2; q/n_1)}{aq^{3/2-3i\tau}}\mathfrak{J}_{2}(q,m,\tau)\Bigr|\mathrm{d}\tau.
\end{align*}
(Recall that $U(x)=1$ for $x\in [1,2]$.) Next we apply Cauchy and Lemma \ref{ram-on-av}, to get
\begin{align}
\label{cauchy}
S_2(N,C)\ll t^\varepsilon N^{1/2}K\int_{-\frac{(NK)^{1/2}}{C}t^\varepsilon}^{+\frac{(NK)^{1/2}}{C}t^\varepsilon}\sum_\pm\sum_{\substack{1\leq L\ll N^{1/2}K^{3/2}t^\varepsilon\\ \text{dyadic}}}\:L^{1/2}\:\left[S_{2,\pm}(N,C,L,\tau)\right]^{1/2}\mathrm{d}\tau
\end{align}
where
\begin{align*}
S_{2,\pm}(N,C,L,\tau)=& \mathop{\sum\sum}_{n_1, n_2} \frac{1}{n_2}U\left(\frac{n_1^2n_2}{L}\right)\Bigl|\mathop{\sum\sum}_{\substack{C< q\leq 2C,\; (m,q)=1\\1\leq |m|\ll qt^{1+\varepsilon}/N\\n_1|q}}\frac{S(\bar m, \pm n_2; q/n_1)}{aq^{3/2-3i\tau}}\mathfrak{J}_{2}(q,m,\tau)\Bigr|^2.
\end{align*}\\

For notational simplicity let us only consider $S_{2,+}(N,C,L,\tau)$. Opening the absolute square and interchanging the order of summations we arrive at following expression for $S_{2,+}(N,C,L,\tau)$ -
\begin{align*}
 \mathop{\sum}_{n_1\leq 2C}\mathop{\sum\sum}_{\substack{C< q\leq 2C,\; (m,q)=1\\1\leq |m|\ll qt^{1+\varepsilon}/N\\n_1|q}}\mathop{\sum\sum}_{\substack{C< q'\leq 2C,\; (m',q')=1\\1\leq |m'|\ll q't^{1+\varepsilon}/N\\n_1|q'}}
\frac{1}{aa'q^{3/2-3i\tau}q'^{3/2+3i\tau}}\mathfrak{J}_{2}(q,m,\tau)\overline{\mathfrak{J}_{2}(q',m',\tau)}\;T.
\end{align*}
where, temporarily,
$$
T=\mathop{\sum}_{n_2} \frac{1}{n_2}U\left(\frac{n_1^2n_2}{L}\right)S(\bar m, n_2; q/n_1)S(\bar m', n_2; q'/n_1).
$$
Let $\hat q=q/n_1$ and $\hat q'=q'/n_1$. Breaking the sum modulo $\hat q\hat q'$ we get
$$
T=\mathop{\sum}_{\beta \bmod{\hat q\hat q'}}S(\bar m, \beta ; \hat q)S(\bar m', \beta; \hat q')\mathop{\sum}_{n_2\in\mathbb Z} \frac{1}{\beta+n_2\hat q\hat q'}U\left(\frac{n_1^2(\beta+n_2\hat q\hat q')}{L}\right).
$$
Applying the Poisson summation formula we have
$$
T=\mathop{\sum}_{\beta \bmod{\hat q\hat q'}}S(\bar m, \beta ; \hat q)S(\bar m', \beta; \hat q')\mathop{\sum}_{n_2\in\mathbb Z} \int_{\mathbb R}\frac{1}{\beta+y\hat q\hat q'}U\left(\frac{n_1^2(\beta+y\hat q\hat q')}{L}\right)e(-n_2y)\mathrm{d}y.
$$
Making the change of variables $n_1^2(\beta+y\hat q\hat q')/L\mapsto w$ it follows that
$$
T=\frac{n_1^2}{qq'}\mathop{\sum}_{n_2\in\mathbb Z} \left[\mathop{\sum}_{\beta \bmod{\hat q\hat q'}}S(\bar m, \beta ; \hat q)S(\bar m', \beta; \hat q')e\left(\frac{\beta n_2}{\hat q\hat q'}\right)\right]\int_{\mathbb R}\frac{1}{w}U\left(w\right)e\left(-\frac{n_2Lw}{qq'}\right)\mathrm{d}w.
$$
The integral is arbitrarily small if $|n_2|\gg C^2t^\varepsilon/L$. (Recall that $a,\: a'\asymp (N/K)^{1/2}$.)

\begin{lemma}
\label{bd-SNCL}
The sum $S_{2,+}(N,C,L,\tau)$ is dominated by the sum
\begin{align*}
\frac{K}{NC^5}B(C,\tau)^2\mathop{\sum}_{n_1\leq 2C}n_1^2\mathop{\sum\sum}_{\substack{C< q\leq 2C,\; (m,q)=1\\1\leq |m|\ll qt^{1+\varepsilon}/N\\n_1|q}}\;\mathop{\sum\sum}_{\substack{C< q'\leq 2C,\; (m',q')=1\\1\leq |m'|\ll q't^{1+\varepsilon}/N\\n_1|q'}}\;\sum_{|n_2|\ll C^2t^\varepsilon/L}|\mathfrak{C}| + O(t^{-2012}),
\end{align*}
where $B(C,\tau)$ is as given in Lemma~\ref{i3-final} (and defined in \eqref{bctau}) and
$$
\mathfrak{C}=\mathop{\sum}_{\beta \bmod{\hat q\hat q'}}S(\bar m, \beta ; \hat q)S(\bar m', \beta; \hat q')e\left(\frac{\beta n_2}{\hat q\hat q'}\right).
$$
\end{lemma}

We have already encountered the character sum $\mathfrak{C}$ in \cite{Mu0}, where we have proved the following. 
\begin{lemma}
\label{char-sum}
We have
$$
\mathfrak{C}\ll \hat q\hat q'(\hat q,\hat q', n_2).
$$
Moreover for $n_2=0$ we get that $\mathfrak{C}=0$ unless $\hat q=\hat q'$, in which case we get
$$
\mathfrak{C}\ll \hat q^2 (\hat q, m-m').
$$
\end{lemma}

For sake of completeness we include the proof here. Let $p$ be a prime, $\hat q=p^jr$ and $\hat q'=p^kr'$ with $p\nmid rr'$. The $p$-part of $\mathfrak{C}$ is given by
$$
\mathfrak{C}_p=\sum_{\beta\bmod{p^{j+k}}}S(\overline{mr}, \beta\bar r; p^j)S(\overline{m'r'}, \beta\bar r'; p^k)e\left(\frac{\beta\overline{rr'}n_2}{p^{j+k}}\right).
$$ 
Opening the Kloosterman sums we get
\begin{align*}
\mathfrak{C}_p=&\mathop{\sideset{}{^\star}\sum\sideset{}{^\star}\sum}_{\substack{\alpha\bmod{p^j}\\\alpha'\bmod{p^k}}}e\left(\frac{\alpha\overline{mr}}{p^j}+\frac{\alpha'\overline{m'r'}}{p^k}\right)\sum_{\beta\bmod{p^{j+k}}}e\left(\frac{\beta\overline{r\alpha}p^k+\beta\overline{r'\alpha'}p^j+\beta\overline{rr'}n_2}{p^{j+k}}\right)\\
=&p^{j+k}\mathop{\sideset{}{^\star}\sum_{\alpha\bmod{p^j}}\;\sideset{}{^\star}\sum_{\alpha'\bmod{p^k}}}_{\overline{r\alpha}p^k+\overline{r'\alpha'}p^j+\overline{rr'}n_2\equiv 0\bmod{p^{j+k}}}e\left(\frac{\alpha\overline{mr}}{p^j}+\frac{\alpha'\overline{m'r'}}{p^k}\right).
\end{align*}
The last sum vanishes unless $(p^j,p^k)|n_2$, and in this case we get that
$$
\mathfrak{C}_p\ll p^{j+k}(p^j,p^k)\ll p^{j+k}(p^j,p^k, n_2).
$$
The first part of the lemma follows. Notice that for $n_2=0$, the congruence has no solutions unless $j=k$, and we conclude the second statement of the lemma.\\


\subsection{Bounding $S_2(N,C)$} Let us first consider the contribution of the zero frequency $n_2=0$ in $S_{2,+}(N,C,L,\tau)$. This we will denote by $S^\flat_{2,+}(N,C,L,\tau)$. Applying the second statement of Lemma~\ref{char-sum}, and Lemma~\ref{bd-SNCL}, we have
\begin{align*}
S^\flat_{2,+}(N,C,L,\tau)\ll \frac{K}{NC^5}B(C,\tau)^2\mathop{\sum}_{n_1\leq 2C}\left\{\frac{C^5t}{n_1^2N}+\frac{C^5t^2}{n_1N^2}\right\}t^{\varepsilon}
\end{align*}
where the first term on the right hand side is the diagonal (i.e. $m=m'$) contribution and the other term is the off-diagonal contribution. Since we are assuming that $N>t$, the diagonal term dominates and we get
\begin{align*}
S^\flat_{2,+}(N,C,L,\tau)\ll \frac{Kt}{N^2}B(C,\tau)^2t^{\varepsilon}.
\end{align*}\\

Next consider the contribution of the non-zero frequencies $n_2\neq 0$ in $S_{2,+}(N,C,L,\tau)$. This we will denote by $S^\sharp_{2,+}(N,C,L,\tau)$. Using the first statement of Lemma~\ref{char-sum} we have
\begin{align*}
S^\sharp_{2,+}(N,C,L,\tau)\ll & \frac{K}{NC^3}B(C,\tau)^2\mathop{\sum}_{n_1\leq 2C}\mathop{\sum\sum}_{\substack{C< q\leq 2C,\; (m,q)=1\\1\leq |m|\ll qt^{1+\varepsilon}/N\\n_1|q}}\;\mathop{\sum\sum}_{\substack{C< q'\leq 2C,\; (m',q')=1\\1\leq |m'|\ll q't^{1+\varepsilon}/N\\n_1|q'}}\;\sum_{1\leq |n_2|\ll C^2t^\varepsilon/L}(q,n_2)\\
\ll & \frac{K}{NC^3}B(C,\tau)^2\mathop{\sum}_{n_1\leq 2C}\left(\frac{C}{n_1}\right)^2\;\left(\frac{Ct}{N}\right)^2\;\frac{C^2}{L}t^\varepsilon \ll \frac{Kt^2C^3}{N^3L}B(C,\tau)^2t^\varepsilon.
\end{align*}
Substituting the above bounds (and the similar bounds for $S_{2,-}(N,C,L,\tau)$) in \eqref{cauchy}, we get
\begin{align*}
S_2(N,C)& \ll t^\varepsilon N^{1/2}K\int_{-\frac{(NK)^{1/2}}{C}t^\varepsilon}^{+\frac{(NK)^{1/2}}{C}t^\varepsilon}B(C,\tau)\sum_\pm\sum_{\substack{1\leq L\ll N^{1/2}K^{3/2}t^\varepsilon\\ \text{dyadic}}}\left\{\frac{(tLK)^{1/2}}{N}+\frac{tC^{3/2}K^{1/2}}{N^{3/2}}\right\}\mathrm{d}\tau\\
&\ll t^\varepsilon \frac{N^{1/2}K^{1/2}}{t^{1/2}}\left\{1+\frac{N}{C^2K^{3/2}}\right\}\left\{\frac{t^{1/2}K^{5/4}}{N^{3/4}}+\frac{tC^{3/2}K^{1/2}}{N^{3/2}}\right\}.
\end{align*}
The second inequality follows from \eqref{bctau-int}. The product of the last two factors is given by
\begin{align*}
\frac{t^{1/2}K^{5/4}}{N^{3/4}}+\frac{tC^{3/2}K^{1/2}}{N^{3/2}}+\frac{t^{1/2}N^{1/4}}{C^2K^{1/4}}+\frac{t}{C^{1/2}N^{1/2}K}.
\end{align*}
Using $N/t^{1+\varepsilon}\ll C\leq (N/K)^{1/2}$, we see that the above sum is dominated by 
\begin{align*}
t^\varepsilon\left(\frac{t^{1/2}K^{5/4}}{N^{3/4}}+\frac{t}{N^{3/4}K^{1/4}}+\frac{t^{5/2}}{N^{7/4}K^{1/4}}+\frac{t^{3/2}}{NK}\right).
\end{align*}
The third term essentially dominates the second as $N\ll t^{3/2+\varepsilon}$. Also the third term dominates the fourth term as $K>N/t^{4/3}$ (which follows from \eqref{cond-K}). Hence we get
\begin{align*}
S_2(N,C)& \ll t^\varepsilon \left\{\frac{K^{7/4}}{N^{1/4}}+\frac{t^{2}K^{1/4}}{N^{5/4}}\right\}.
\end{align*}
Substituting this estimate in Lemma~\ref{pre-cauchy} and Lemma~\ref{pre-lemma}, we get that the contribution of $S_2(N,C)$ to $S^+(N)$ is bounded by
\begin{align}
\label{s2nc-cont}
t^\varepsilon N^{3/4}\left\{K^{3/4}+\frac{t^2}{NK^{3/4}}\right\}.
\end{align}


\section{Application of Cauchy inequality and Poisson summation - II}

It remains to estimate $S_{1,J}(N.C)$. This is comparatively delicate as we need to get cancellation in the integral over $\tau$ for large $J$. For notational simplicity let us only consider the case of positive $J$ with $J\gg t^\varepsilon$. The same analysis holds for negative  $J$ with $-J\gg t^\varepsilon$. For $J$ of smaller size, the analysis is even simpler as there is no need to get cancellation in the $\tau$ integral.

\subsection{Applying Cauchy inequality and Poisson summation} As before we take dyadic segmentation, but keep the integral over $\tau$ inside the absolute value to get
\begin{align*}
S_{1,J}&(N,C)\leq N^{1/2}K\sum_\pm \sum_{\substack{1\leq L\ll N^{1/2}K^{3/2}t^\varepsilon\\ \text{dyadic}}} \mathop{\sum\sum}_{n_1, n_2} \frac{|\lambda(n_2,n_1)|}{n_2^{1/2}}U\left(\frac{n_1^2n_2}{L}\right)\\
&\times \Bigl|\int_{\mathbb R}\left(n_1^2n_2N\right)^{-i\tau}\gamma_\pm\left(-\frac{1}{2}+i\tau\right)\mathop{\sum\sum}_{\substack{C< q\leq 2C,\; (m,q)=1\\1\leq |m|\ll qt^{1+\varepsilon}/N\\n_1|q}}\frac{S(\bar m, \pm n_2; q/n_1)}{aq^{3/2-3i\tau}}\mathfrak{J}_{1}(q,m,\tau)W_J(\tau)\mathrm{d}\tau\Bigr|.
\end{align*} 
Applying Cauchy and Lemma \ref{ram-on-av}, we conclude that
\begin{align}
\label{cauchy-2}
S_{1,J}(N,C)\ll t^\varepsilon N^{1/2}K\sum_\pm\sum_{\substack{1\leq  L\ll N^{1/2}K^{3/2}t^\varepsilon\\ \text{dyadic}}}\:L^{1/2}\:\left[S_{1,J,\pm}(N,C,L)\right]^{1/2}
\end{align}
where $S_{1,J,\pm}(N,C,L)$ is given by
\begin{align*}
\mathop{\sum\sum}_{n_1, n_2} \frac{1}{n_2}U\left(\frac{n_1^2n_2}{L}\right)\Bigl|\int_{\mathbb R}&\left(n_1^2n_2N\right)^{-i\tau}\gamma_\pm\left(-\frac{1}{2}+i\tau\right)\\
&\mathop{\sum\sum}_{\substack{C< q\leq 2C,\; (m,q)=1\\1\leq |m|\ll qt^{1+\varepsilon}/N\\n_1|q}}\frac{S(\bar m, \pm n_2; q/n_1)}{aq^{3/2-3i\tau}}\mathfrak{J}_{1}(q,m,\tau)W_J(\tau)\mathrm{d}\tau\Bigr|^2.
\end{align*}\\

We will only consider $S_{1,J,+}(N,C,L)$. Opening the absolute square and interchanging the order of summations we arrive at
\begin{align*}
S_{1,J,+}&(N,C,L)= \mathop{\sum}_{n_1\leq 2C}\mathop{\iint}_{\mathbb R^2}\left(n_1^2N\right)^{-i(\tau-\tau')}\gamma_+\left(-\frac{1}{2}+i\tau\right)\overline{\gamma_+\left(-\frac{1}{2}+i\tau'\right)}W_J(\tau)W_J(\tau')\\
&\times \mathop{\sum\sum}_{\substack{C< q\leq 2C,\; (m,q)=1\\1\leq |m|\ll qt^{1+\varepsilon}/N\\n_1|q}}\mathop{\sum\sum}_{\substack{C< q'\leq 2C,\; (m',q')=1\\1\leq |m'|\ll q't^{1+\varepsilon}/N\\n_1|q'}}
\frac{1}{aa'q^{\frac{3}{2}-3i\tau}q'^{\frac{3}{2}+3i\tau'}}\mathfrak{J}_{1}(q,m,\tau)\overline{\mathfrak{J}_{1}(q',m',\tau')}\;T\;\mathrm{d}\tau\mathrm{d}\tau'.
\end{align*}
where (temporarily)
$$
T=\mathop{\sum}_{n_2\in \mathbb Z} n_2^{-1-i(\tau-\tau')}U\left(\frac{n_1^2n_2}{L}\right)S(\bar m, n_2; q/n_1)S(\bar m', n_2; q'/n_1).
$$
Let $\hat q=q/n_1$ and $\hat q'=q'/n_1$. Breaking the sum modulo $\hat q\hat q'$ and applying Poisson summation we get
$$
T=\mathop{\sum}_{\beta \bmod{\hat q\hat q'}}S(\bar m, \beta ; \hat q)S(\bar m', \beta; \hat q')\mathop{\sum}_{n_2\in\mathbb Z} \int_{\mathbb R}(\beta+y\hat q\hat q')^{-1-i(\tau-\tau')}U\left(\frac{n_1^2(\beta+y\hat q\hat q')}{L}\right)e(-n_2y)\mathrm{d}y.
$$
Making the change of variables $n_1^2(\beta+y\hat q\hat q')/L\mapsto w$ it follows that
$$
T=\frac{n_1^2}{qq'}\left(\frac{L}{n_1^2}\right)^{-i(\tau-\tau')}\mathop{\sum}_{n_2\in\mathbb Z} \mathfrak{C}\;U^{ \dagger}\left(\frac{n_2L}{qq'},-i(\tau-\tau')\right),
$$
where $\mathfrak{C}$ is the same character sum that appears in Lemma~\ref{bd-SNCL}, and the exponential integral $U^{ \dagger}$ is as defined in Subsection~\ref{anintegral}. From the second statement of Lemma~\ref{the-integral} we see that the integral is arbitrarily small if $|n_2|\gg C(NK)^{1/2}t^\varepsilon/L$. (Recall that $|\tau-\tau'|\ll (NK)^{1/2}t^\varepsilon/C$ and $q, q'\sim C$.)\\ 

\begin{lemma}
\label{lemma-for-s1nc}
The sum $S_{1,J,+}(N,C,L)$ is dominated by the sum
\begin{align*}
\frac{K}{NC^5}\mathop{\sum}_{n_1\leq 2C}n_1^2\mathop{\sum\sum}_{\substack{C< q\leq 2C,\; (m,q)=1\\1\leq |m|\ll qt^{1+\varepsilon}/N\\n_1|q}}\;\mathop{\sum\sum}_{\substack{C< q'\leq 2C,\; (m',q')=1\\1\leq |m'|\ll q't^{1+\varepsilon}/N\\n_1|q'}}\;\sum_{|n_2|\ll C(NK)^{1/2}t^\varepsilon/L}|\mathfrak{C}||\mathfrak{K}|+O(t^{-2012}),
\end{align*}
where $\mathfrak{C}$ is as in Lemma~\ref{bd-SNCL}, and 
\begin{align*}
\mathfrak{K}=\mathop{\iint}_{\mathbb R^2}&\gamma_+\left(-\frac{1}{2}+i\tau\right)\overline{\gamma_+\left(-\frac{1}{2}+i\tau'\right)}\frac{\left(LN\right)^{-i(\tau-\tau')}}{q^{-3i\tau}q'^{3i\tau'}}W_J(\tau)W_J(\tau')\\
\times &\mathfrak{J}_{1}(q,m,\tau)\overline{\mathfrak{J}_{1}(q',m',\tau')}\:U^{ \dagger}\left(\frac{n_2L}{qq'},-i(\tau-\tau')\right)\;\mathrm{d}\tau\mathrm{d}\tau'.
\end{align*}\\
\end{lemma}

We already have a satisfactory bound for the character sum $\mathfrak{C}$. We only need to estimate the exponential integral $\mathfrak{K}$. Using the explicit form of $\mathfrak{J}_{1}(q,m,\tau)$, as given in Lemma~\ref{i3-final}, we get
\begin{align*}
\mathfrak{K}=&\frac{|c_4|^2}{K^2}\mathop{\iint}_{\mathbb R^2}\gamma_+\left(-\frac{1}{2}+i\tau\right)\overline{\gamma_+\left(-\frac{1}{2}+i\tau'\right)}\:W_J(q,m,\tau)W_J(q',m',\tau')\frac{\left(LN\right)^{-i(\tau-\tau')}}{q^{-3i\tau}q'^{3i\tau'}}\\
\times &\left(-\frac{(t+\tau)q}{2\pi eNm}\right)^{-i(t+\tau)}\left(-\frac{(t+\tau')q'}{2\pi eNm'}\right)^{i(t+\tau')}U^{ \dagger}\left(\frac{n_2L}{qq'},-i(\tau-\tau')\right)\mathrm{d}\tau\mathrm{d}\tau',
\end{align*}
where
$$
W_J(q,m,\tau)=(t+\tau)^{-1/2}W_J(\tau)\left(-\frac{(t+\tau)q}{2\pi eNm}\right)^{3/2}V\left(-\frac{(t+\tau)q}{2\pi Nm}\right)\:\int_0^1V\left(\frac{\tau}{K}-\frac{(t+\tau)x}{Kma}\right)\mathrm{d}x.
$$
Since $u^{3/2}V(u)\ll 1$, and $|\tau|\ll J\ll t^{1-\varepsilon}$, it follows that 
$$
\frac{\partial}{\partial \tau}W_J(q,m,\tau)\ll \frac{1}{t^{1/2}|\tau|}.
$$\\


\subsection*{The integral $\mathfrak{K}$} Let us first consider the integral 
\begin{align*}
U^{\dagger}\left(\frac{n_2L}{qq'},-i(\tau-\tau')\right)=\int_{\mathbb R}U\left(w\right)w^{-i(\tau-\tau')-1}e\left(-\frac{n_2Lw}{qq'}\right)\mathrm{d}w
\end{align*}
which appears in Lemma~\ref{lemma-for-s1nc}. We study this in the light of Lemma~\ref{the-integral}. For $n_2=0$ the integral is negligibly small if $|\tau-\tau'|\gg t^\varepsilon$. So in this case we get 
$$
\mathfrak{K}\ll \frac{N^{1/2}t^\varepsilon}{K^{3/2}Ct}.
$$ \\

Now suppose $n_2\neq 0$. In this case we apply the first statement of Lemma~\ref{the-integral} to deduce 
\begin{align*}
U^{\dagger}\left(\frac{n_2L}{qq'},-i(\tau-\tau')\right)=\frac{c_5}{(\tau'-\tau)^{1/2}}&U\left(\frac{(\tau'-\tau)qq'}{2\pi n_2L}\right)\left(\frac{(\tau'-\tau)qq'}{2\pi e n_2L}\right)^{-i(\tau-\tau')}\\
&+O\left(\min\left\{\frac{1}{|\tau-\tau'|^{3/2}},\frac{C^3}{(|n_2|L)^{3/2}}\right\}\right)
\end{align*}
for some constant $c_5$ (which depends on the sign of $n_2$). The contribution of the error term towards $\mathfrak{K}$ is bounded by
$$
O\left(\frac{1}{K^2 t}\mathop{\iint}_J^{2J} \min\left\{\frac{1}{|\tau-\tau'|^{3/2}},\frac{C^3}{(|n_2|L)^{3/2}}\right\}\mathrm{d}\tau\mathrm{d}\tau'\right).
$$
Now we see that
$$
\frac{1}{K^2 t}\mathop{\iint}_{\substack{[J,2J]^2\\|\tau-\tau'|\leq |n_2L|/C^2}} \frac{C^3}{(|n_2|L)^{3/2}}\mathrm{d}\tau\mathrm{d}\tau'\ll \frac{C}{K^2 t(|n_2|L)^{1/2}}J\ll \frac{1}{K^{3/2}t}\frac{N^{1/2}}{(|n_2|L)^{1/2}}t^\varepsilon,
$$
(as $J\ll (NK)^{1/2}t^\varepsilon/C$) and
\begin{align*}
\frac{1}{K^2 t}\mathop{\iint}_{\substack{[J,2J]^2\\|\tau-\tau'|> |n_2L|/C^2}} &\frac{1}{|\tau-\tau'|^{3/2}}\mathrm{d}\tau\mathrm{d}\tau'\ll \frac{1}{K^2 t}\frac{C}{(|n_2|L)^{1/2}}\mathop{\iint}_{\substack{[J,2J]^2}} \frac{1}{|\tau-\tau'|^{1-\varepsilon}}\mathrm{d}\tau\mathrm{d}\tau'\\
&\ll \frac{C}{K^2 t(|n_2|L)^{1/2}}Jt^\varepsilon\ll \frac{1}{K^{3/2}t}\frac{N^{1/2}}{(|n_2|L)^{1/2}}t^\varepsilon.
\end{align*}
We set 
$$
B^\star(C,0)=\frac{N^{1/2}}{K^{3/2}Ct}
$$ 
and for $n_2\neq 0$
\begin{align}
\label{b'tau}
B^\star(C,n_2)=\frac{1}{K^{3/2}t}\frac{N^{1/2}}{(|n_2|L)^{1/2}}.
\end{align}
Then for $q\neq 0$ (and $C\leq (N/K)^{1/2}$) we have
\begin{align}
\label{b'tau-sum}
\sum_{1\leq |n_2|\ll C(NK)^{1/2}t^\varepsilon/L}(q,n_2)B^\star(C,n_2) \ll\frac{N}{K^{3/2}tL}t^\varepsilon.
\end{align}\\

Now we consider the main term. We pull out the oscillation from the gamma factors using \eqref{scaling}. By Fourier inversion we write
$$
\left(\frac{2\pi n_2L}{(\tau'-\tau)qq'}\right)^{1/2}U\left(\frac{(\tau'-\tau)qq'}{2\pi n_2L}\right)=\int_{\mathbb R}U^{\dagger}\left(r, \frac{1}{2}\right)e\left(\frac{(\tau'-\tau)qq'}{2\pi n_2L}r\right)\mathrm{d}r.
$$ 
We conclude that (for some constant $c_6$ depending on the sign of $n_2$)
\begin{align}
\label{lead-term-for-K}
&\mathfrak{K}=\frac{c_6}{K^2}\left(\frac{qq'}{|n_2|L}\right)^{1/2}\int_{\mathbb R}U^{\dagger}\left(r, \frac{1}{2}\right)\mathop{\iint}_{\mathbb R^2}g(\tau,\tau')e(f(\tau,\tau'))\mathrm{d}\tau\mathrm{d}\tau'\mathrm{d}r + O(B^\star(C,n_2)t^\varepsilon).
\end{align}
where
\begin{align*}
2\pi f(\tau,\tau')=3\tau&\log\left(\frac{\tau}{e\pi} \right)-3\tau'\log\left(\frac{\tau'}{e\pi} \right)-(\tau-\tau')\log LN+3\tau\log q-3\tau'\log q'\\
&-(t+\tau)\log\left(-\frac{(t+\tau)q}{2\pi eNm}\right)+(t+\tau')\log\left(-\frac{(t+\tau')q'}{2\pi eNm'}\right)\\
&-(\tau-\tau')\log\left(\frac{(\tau'-\tau)qq'}{2\pi e n_2L}\right)+\frac{(\tau'-\tau)qq'}{n_2L}r.
\end{align*}
and 
\begin{align*}
&g(\tau,\tau')=\Phi_+\left(\tau\right)\overline{\Phi_+\left(\tau'\right)}W_J(q,m,\tau)W_J(q',m',\tau').
\end{align*}
We will use Lemma~\ref{exp-2-var-lemma} to analyse the double exponential integral over $\tau$ and $\tau'$. Differentiating we get
\begin{align*}
2\pi \frac{\partial^2}{\partial \tau^2}f(\tau,\tau')=&\frac{3}{\tau}-\frac{1}{t+\tau}+\frac{1}{\tau'-\tau},\;\;\;\;\;2\pi \frac{\partial^2}{\partial \tau^{'2}}f(\tau,\tau')=-\frac{3}{\tau'}+\frac{1}{t+\tau'}+\frac{1}{\tau'-\tau}
\end{align*}
and
\begin{align*}
2\pi \frac{\partial^2}{\partial \tau' \partial \tau}f(\tau,\tau')=-\frac{1}{\tau'-\tau}.
\end{align*}
Also by explicit computation we get (using that $J\ll t^{1-\varepsilon}$)
\begin{align*}
4\pi^2\left[\frac{\partial^2}{\partial \tau^2}f(\tau,\tau')\frac{\partial^2}{\partial \tau^{'2}}f(\tau,\tau')-\left(\frac{\partial^2}{\partial \tau' \partial \tau}f(\tau,\tau')\right)^2\right]=-\frac{6}{\tau\tau'}+O\left(\frac{1}{tJ}\right),
\end{align*}
for $\tau$, $\tau'$ such that $g(\tau,\tau')\neq 0$. (Recall that $\mathrm{Supp}W_J\subset [J,4J/3]$.) So the conditions of the Lemma~\ref{exp-2-var-lemma} hold with $r_1=r_2=1/J^{1/2}$. Next we need to compute the total variation of the weight function $g(\tau,\tau')$. Recall that $\Phi_+'(\tau)\ll |\tau|^{-1}$ and $W_J'(q,m,\tau)\ll t^{-1/2}|\tau|^{-1}$ (derivative with respect to $\tau$). It follows that
$
\mathrm{var}(g)\ll t^{-1+\varepsilon}.
$
So from Lemma~\ref{exp-2-var-lemma} we conclude that the double integral (over $\tau$, $\tau'$) is bounded by
$
O\left(Jt^{-1+\varepsilon}\right).
$
Then integrating trivially over $r$ using the rapid decay of the Fourier transform we get that the total contribution of the leading term in \eqref{lead-term-for-K} to $\mathfrak{K}$ is bounded by
$$
O\left(\frac{1}{K^2}\frac{C}{(|n_2|L)^{1/2}}\frac{(NK)^{1/2}}{C}t^{-1+\varepsilon}\right)=O(B^\star(C,n_2) t^\varepsilon).
$$\\

\begin{lemma}
\label{bd-fin}
We have
$$
\mathfrak{K}\ll B^\star(C,n_2)t^\varepsilon
$$
where $B^\star(C,n_2)$ is given by \eqref{b'tau}.
\end{lemma}


\subsection{Bounding $S_{1,J,\pm}(N,C,L)$ and $S_{1,J}(N,C)$} From Lemma~\ref{lemma-for-s1nc} and Lemma~\ref{bd-fin}, it follows that we need to estimate the sum 
\begin{align*}
\frac{K}{NC^5}\mathop{\sum}_{n_1\leq 2C}n_1^2\mathop{\sum\sum}_{\substack{C< q\leq 2C,\; (m,q)=1\\1\leq |m|\ll qt^{1+\varepsilon}/N\\n_1|q}}\;\mathop{\sum\sum}_{\substack{C< q'\leq 2C,\; (m',q')=1\\1\leq |m'|\ll q't^{1+\varepsilon}/N\\n_1|q'}}\;\sum_{|n_2|\ll C(NK)^{1/2}t^\varepsilon/L}|\mathfrak{C}|B^\star(C,n_2).
\end{align*}
Let us first consider the contribution of the zero frequency $n_2=0$. We denote this by $S^\flat_{1,J,+}(N,C,L)$. Applying the second statement of Lemma~\ref{char-sum}, we have
\begin{align*}
S^\flat_{1,J,+}(N,C,L)\ll \frac{K}{NC^5}\frac{N^{1/2}}{K^{3/2}Ct}\mathop{\sum}_{n_1\leq 2C}\left\{\frac{C^5t}{n_1^2N}+\frac{C^5t^2}{n_1N^2}\right\}t^{\varepsilon}
\end{align*}
where the first term on the right hand side is the diagonal (i.e. $m=m'$) contribution and the other term is the off-diagonal contribution. Since we are assuming that $N>t$, the diagonal term dominates and we get
\begin{align*}
S^\flat_{1,J,+}(N,C,L)\ll \frac{1}{N^{3/2}K^{1/2}C}t^{\varepsilon}\ll \frac{t}{N^{5/2}K^{1/2}}t^{\varepsilon}.
\end{align*}
(Recall that $C>N/t^{1+\varepsilon}$.)\\

Next consider the contribution of the non-zero frequencies $n_2\neq 0$ in $S_{1,J,+}(N,C,L)$. This we will denote by $S^\sharp_{1,J,+}(N,C,L)$. We have, using \eqref{b'tau-sum},
\begin{align*}
S^\sharp_{1,J,+}(N,C,L)\ll & \frac{Kt^\varepsilon}{NC^3}\mathop{\sum}_{n_1\leq 2C}\mathop{\sum\sum}_{\substack{C< q\leq 2C,\; (m,q)=1\\1\leq |m|\ll qt^{1+\varepsilon}/N\\n_1|q}}\;\mathop{\sum\sum}_{\substack{C< q'\leq 2C,\; (m',q')=1\\1\leq |m'|\ll q't^{1+\varepsilon}/N\\n_1|q'}}\;\sum_{1\leq |n_2|\ll C(NK)^{1/2}t^\varepsilon/L}(q,n_2)B^\star(C,n_2)\\
\ll & \frac{Kt^\varepsilon}{NC^3}\mathop{\sum}_{n_1\leq 2C}\left(\frac{C}{n_1}\right)^2\;\left(\frac{Ct}{N}\right)^2\;\frac{N}{K^{3/2}tL}\ll \frac{tC}{N^2K^{1/2}L}t^\varepsilon.
\end{align*}
We conclude that
\begin{align*}
S_{1,J,+}(N,C,L)\ll \left(\frac{t}{N^{5/2}K^{1/2}}+\frac{tC}{N^2K^{1/2}L}\right)t^\varepsilon.
\end{align*}
The same bound holds for $S_{1,J,-}(N,C,L)$.\\

Substituting the above bounds in \eqref{cauchy-2} we get
\begin{align*}
S_{1,J}(N,C)& \ll t^\varepsilon N^{1/2}K\sum_\pm\sum_{\substack{1\leq L\ll N^{1/2}K^{3/2}t^\varepsilon\\ \text{dyadic}}}\left\{\frac{(tL)^{1/2}}{N^{5/4}K^{1/4}}+\frac{(tC)^{1/2}}{NK^{1/4}}\right\}\\
&\ll t^\varepsilon \left\{\frac{K^{3/2}t^{1/2}}{N^{1/2}}+\frac{K^{1/2}t^{1/2}}{N^{1/4}}\right\}.
\end{align*}
The same bound holds for all values of $J$. Since there are $O(\log t)$ many $J$ we can sum over them without worsening the bound, and so the same bound holds for $S_1(N,C):=\sum_JS_{1,J}(N,C)$. Substituting this estimate in Lemma~\ref{pre-lemma} we get that the contribution of $S_1(N,C)$ to $S^+(N)$ is bounded by
\begin{align}
\label{s1nc-cont}
t^\varepsilon N^{3/4}\left\{\frac{K^{1/2}t^{1/2}}{N^{1/4}}+\frac{t^{1/2}}{K^{1/2}}\right\}.
\end{align}\\

Next we combine the bounds from \eqref{s2nc-cont} and \eqref{s1nc-cont} to get a bound for $S^+(N)$. But observe that the first term of \eqref{s1nc-cont} dominates the first term of \eqref{s2nc-cont} as $K<t^2/N$, and it dominates the second term of  \eqref{s2nc-cont} as $K>(t^2/N)^{3/5}$ (see \eqref{cond-K}). So we conclude that
\begin{align*}
S^+(N)\ll t^\varepsilon N^{3/4}t^{1/2}\left\{\frac{K^{1/2}}{N^{1/4}}+\frac{1}{K^{1/2}}\right\}.
\end{align*}
The optimal choice for $K$ is obtained by equating the terms inside the braces. This is given by  
$
K=N^{1/4}.
$
We see that this choice satisfies the imposed condition \eqref{cond-K} on $K$  if $N>t^{24/17}$. In the range $t^{11/8}<N<t^{24/17}$ we pick $K=t^{6/5}/N^{3/5}$. The above bound boils down to
\begin{align*}
S^+(N)\ll \begin{cases} t^{1/2+\varepsilon} N^{5/8} &\text{if}\;\;\;t^{24/17}<N\ll t^{3/2+\varepsilon};\\
t^{11/10+\varepsilon}N^{1/5} & \text{if}\;\;\;t^{11/8}<N\leq t^{24/17}.\end{cases}
\end{align*}
This completes the proof of Proposition~\ref{prop1}.
 

\end{document}